\documentclass[12pt]{article}
\pdfoutput=1
\usepackage{graphicx}
\usepackage{cite}
\graphicspath{{figures/Appendix/}}
\usepackage{amssymb} 
\usepackage{amsmath} 

\begin{document}
\thispagestyle{plain}
\pagestyle{plain}
\pagenumbering{arabic}
\setcounter{page}{1}
\setcounter{section}{0}

\begin{center}
{\bf Projections Onto Convex Sets (POCS) Based Optimization by Lifting}
\end{center}
\begin{center}
A. Enis Cetin, Alican Bozkurt, Osman Gunay, Y. Hakan Habiboglu, Kivanc Kose, Ibrahim Onaran, R. A. Sevimli\\
Dept. of Electrical and Electronic Engineering\\
Bilkent University,\\
 Ankara, Turkey\\
E-mail: cetin at bilkent.edu.tr
\end{center}

\label{app:pocs} 
{\bf Abstract}:

Two new optimization techniques based on projections onto convex space (POCS) framework for solving convex and some non-convex optimization problems are presented. The dimension of the minimization problem is lifted by one and  sets corresponding to the cost function are defined. If the cost function is a convex function in $R^N$ the corresponding set is a convex set in $R^{N+1}$. The iterative optimization approach starts with an arbitrary initial estimate in $R^{N+1}$ and an orthogonal projection is performed onto one of the sets in a sequential manner at each step of the optimization problem. The method provides globally optimal solutions in total-variation, filtered variation, $l_1$, and entropic cost functions. It is also experimentally observed that cost functions based on  $l_p,$ $p<1$ can be handled by using the supporting hyperplane concept.   

\section{Introduction}
In many inverse signal and image processing problems and compressing sensing problems an optimization problem is solved to find a solution:
\begin{equation}
\label{app:eq:c1}
\underset{\mathbf{w}\in \text{C}}{\text{min}} f(\mathbf{w})
\end{equation}
where $C$ is a set in $\mathbb{R}^N$ and $f(\mathbf{w})$ is the cost function. 
Some commonly used cost functions are based on $l_1$, $l_2$, total-variation, filtered variation, and entropic functions \cite{Rud92, Bar07,Can08, Kos12,Gunay}.
%
Bregman developed iterative methods based on the so-called Bregman distance  to solve convex optimization problems which arise in signal and image processing \cite{Bre67}. In Bregman's approach, it is necessary to perform a D-projection (or Bregman projection) at each step of the algorithm an it may not be easy to compute the Bregman distance in general \cite{Yin08,Kiv12,Gunay}.

In this article Bregman's projections onto convex sets (POCS) framework \cite{Bregman,You82} is used to solve convex and some non-convex optimization problems instead of his Bregman distance approach. Bregman's POCS method is widely used for finding a common point of convex sets in many inverse signal and image processing problems\cite{You82,Her95,Cen12,Sla08,Cet03,Cetin94,Cetin89,Kose11,Cen81,Sla09,The11,censor1987optimization,Tru85,Com04,Com93,Kim92,yamada2011minimizing,censor1987some,Sez82,censor1992proximal,Tuy81,censor1981row,censor1991optimization,Ros13}. In the ordinary POCS approach the goal is simply to find a vector which is in the intersection of convex sets.  In each step of the iterative algorithm an orthogonal projection is performed onto one of the convex sets. Bregman showed that successive orthogonal projections converge to a vector which is in the intersection of all the convex sets. If the sets do not intersect iterates oscillate between members of the sets \cite{Gub67,Com12,Cet97}. Since  there is no need to compute the Bregman distance in standard POCS, it found applications in many practical problems.

In our approach the dimension of the minimization problem is lifted by one and  sets corresponding to the cost function are defined. This approach is graphically illustrated in Figure 1. If the cost function is a convex function in $R^N$ the corresponding set is a convex set in $R^{N+1}$. As a result the convex minimization problem is reduced to finding a specific member (the optimal solution) of the set corresponding to the cost function. As in ordinary POCS approach 
the new iterative optimization method starts with an arbitrary initial estimate in $R^{N+1}$ and an orthogonal projection is performed onto one of the sets. After this vector is calculated it is projected onto the other set. This process is continued in a sequential manner at each step of the optimization problem. 
The method provides globally optimal solutions in total-variation, filtered variation, $l_1$, and entropic function based cost functions because they are convex cost functions. It is also experimentally observed that cost functions based on  $l_p,$ $p<1$ can be handled by using the supporting hyperplane concept.   

The article is organized as follows. In Section 2, the convex minimization method based on the POCS approach is introduced. In Section 3, another convex minimization method based on supporting hyperplanes is studied. Since it is very easy to perform an orthogonal projection onto a hyperplane this method is computationally implementable for many cost functions without solving any nonlinear equations. In Section 4, we present some examples on non-convex minimization.

\section{Convex Minimization}

Let us first consider a convex minimization problem 
\begin{equation}
\label{app:eq:c5}
\underset{\mathbf{w}\in \mathbb{R}^N}{\text{min}} f(\mathbf{w})
\end{equation}
where $f:\mathbb{R}^N \rightarrow \mathbb{R}$ is a convex function.

We increase the dimension by one to define the following sets in $\mathbb{R}^{N+1}$ corresponding to the cost function $f(\mathbf{w})$ as follows:
\begin{equation}
\label{app:eq:c6}
\text{C}_f = \{\underline{\mathbf{w}} = [\mathbf{w}^T~y]^T : \mathrm{~} y\geq f(\mathbf{w})\}
\end{equation}
which is the set of $N+1$ dimensional vectors whose $N+1^{st}$ component $y$ is greater than $f(\mathbf{w})$. We use bold face letters for $N$ dimensional vectors and underlined bold face letters for $N+1$ dimensional vectors, respectively.

\begin{figure}[ht!]
\includegraphics[width=\linewidth,clip,trim=0.25in 0.25 0.25 0.25in]{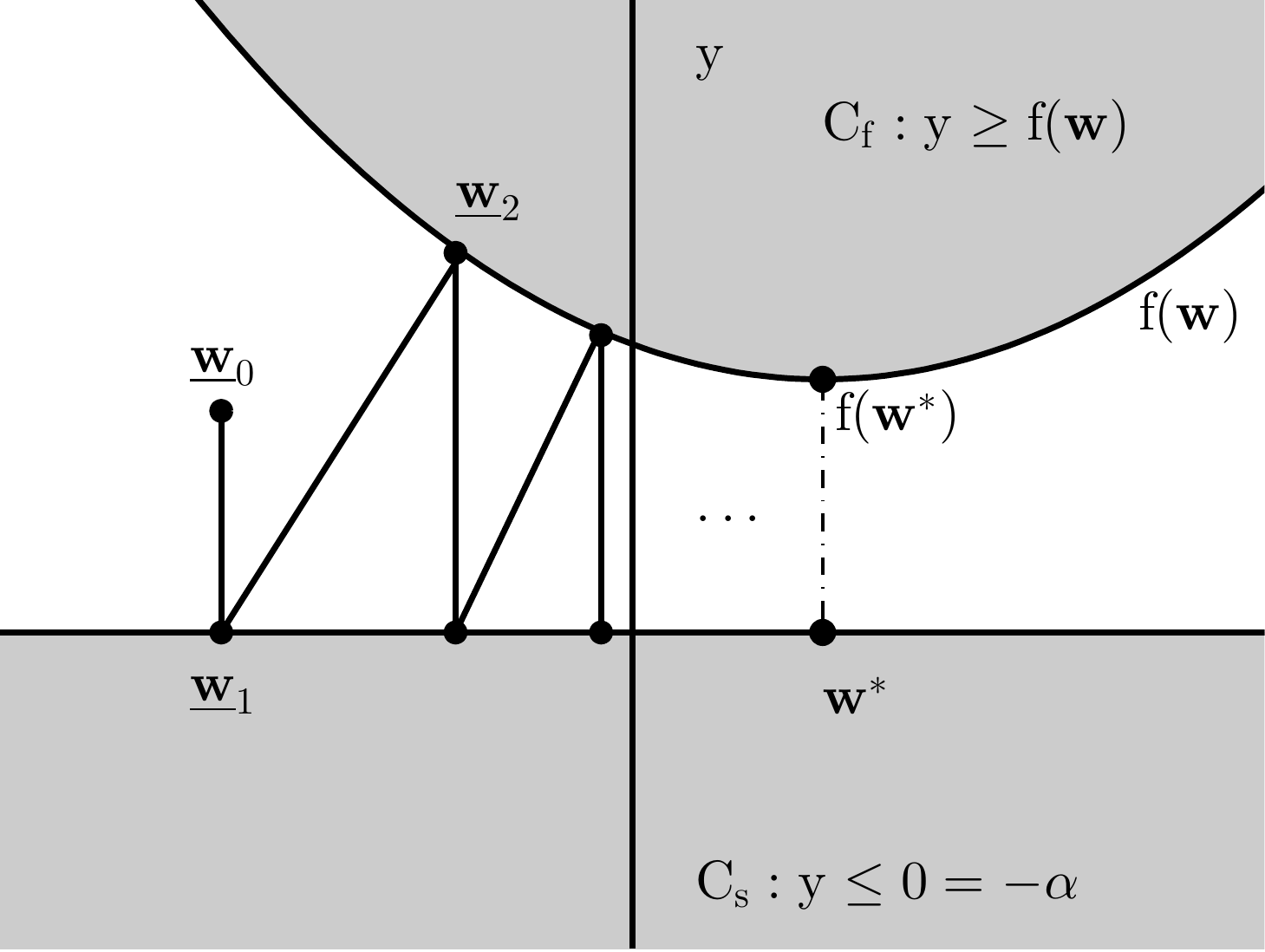}
\caption[Two projecting convex sets.]{Two convex sets C$_{\mathrm{f}}$ and C$_{\mathrm{s}}$ corresponding to the cost function $f$. We sequentially project an initial vector $\underline{\mathbf{w}}_0$ onto $C_s$ and $C_f$ to find the global minimum which is located at $\mathbf{w}^*$.}
\label{app:convex}
\end{figure}

The second set that is related with the cost function $f(\mathbf{w})$ is the level set:
\begin{equation}
\label{app:eq:c7}
\text{C}_s = \{\underline{\mathbf{w}} = [\mathbf{w}^T~y]^T : ~ y\leq \alpha , ~\underline{\mathbf{w}} \in \mathbb{R}^{N+1}\}
\end{equation}
where $\alpha$ is a real number. Here it is assumed that $f(\mathbf{w})\geq \alpha$ for all $f(\mathbf{w})\in \mathbb{R}$ such that the sets C$_{\mathrm{f}}$ and C$_{\mathrm{s}}$ do not intersect. They are both closed and convex sets in  $\mathbb{R}^{N+1}$. Sets C$_{\mathrm{f}}$ and C$_{\mathrm{s}}$ are graphically illustrated in Fig. \ref{app:convex} in which $\alpha=0.$

The POCS based minimization algorithm starts with an arbitrary $\underline{\mathbf{w}}_0 =[ \mathbf{w}_0^T ~ y_0]^T \in \mathbb{R}^{N+1}$. We project $\underline{\mathbf{w}}_0$ onto the set C$_{\mathrm{s}}$ to obtain the first iterate $\underline{\mathbf{w}}_1$ which will be,
\begin{equation}
\underline{\mathbf{w}}_1 = [~\mathbf{w}_0^T~~~0~]^T
\end{equation}
where $\alpha=0$ is assumed as in Fig. \ref{app:convex}. Then we project $\underline{\mathbf{w}}_1$ onto the set C$_{\mathrm{f}}$. The new iterate  $\underline{\mathbf{w}}_2$ is determined by minimizing the distance between $\underline{\mathbf{w}}_1$ and C$_{\mathrm{f}}$, i.e.,

\begin{equation}
\label{app:eq:convex}
\underline{\mathbf{w}}_2 = \text{arg} \underset{\underline{\mathbf{w}}\in \text{C}_{\mathrm{s}}}{\text{min}} \|\underline{\mathbf{w}}_1 - \underline{\mathbf{w}}\|
\end{equation}
Eq. \ref{app:eq:convex} is the ordinary orthogonal projection operation onto the set $\mathrm{C_f} \in \mathbb{R}^{N+1}$. 
To solve the problem in Eq. \ref{app:eq:convex} we do not need to compute the Bregman's so-called D-projection. After finding $\underline{\mathbf{w}}_2$, we perform the next projection onto the set C$_{\mathrm{s}}$ and obtain $\underline{\mathbf{w}}_3$ etc. Eventually iterates oscillate between two nearest vectors of the two sets C$_{\mathrm{s}}$ and C$_{\mathrm{f}}$. As a result we obtain
\begin{equation}
\label{app:eq:convex2}
\underset{n \rightarrow \infty}{\text{lim}} \underline{\mathbf{w}}_{2n} = [~\mathbf{w}^*~~~f(\mathbf{w}^*)~]^T
\end{equation}
where $\mathbf{w}^*$ is the N dimensional vector minimizing $f(\mathbf{w})$. The proof of Equation (\ref{app:eq:convex2}) follows from Bregman's POCS theorem \cite{Bregman,Gub67}. It was generalized to non-intersection case by Gubin et. al \cite{Gub67,Cen12},\cite{Com12}. Since the two closed and convex sets $C_s$ and $C_f$ are closest to each other at the optimal solution case, iterates oscilate between
the vectors $[~\mathbf{w}^*~~~f(\mathbf{w}^*)~]^T$ and $[~\mathbf{w}^*~~~0~]^T$ in $R^{N+1}$ as $n$ tends to infinity. It is possible to increase the speed of convergence by non-orthogonal projections \cite{Com93}.
 
\begin{figure}[ht!]
\includegraphics[width=\linewidth,clip,trim=0.25in 0.25 0.25 0.25in]{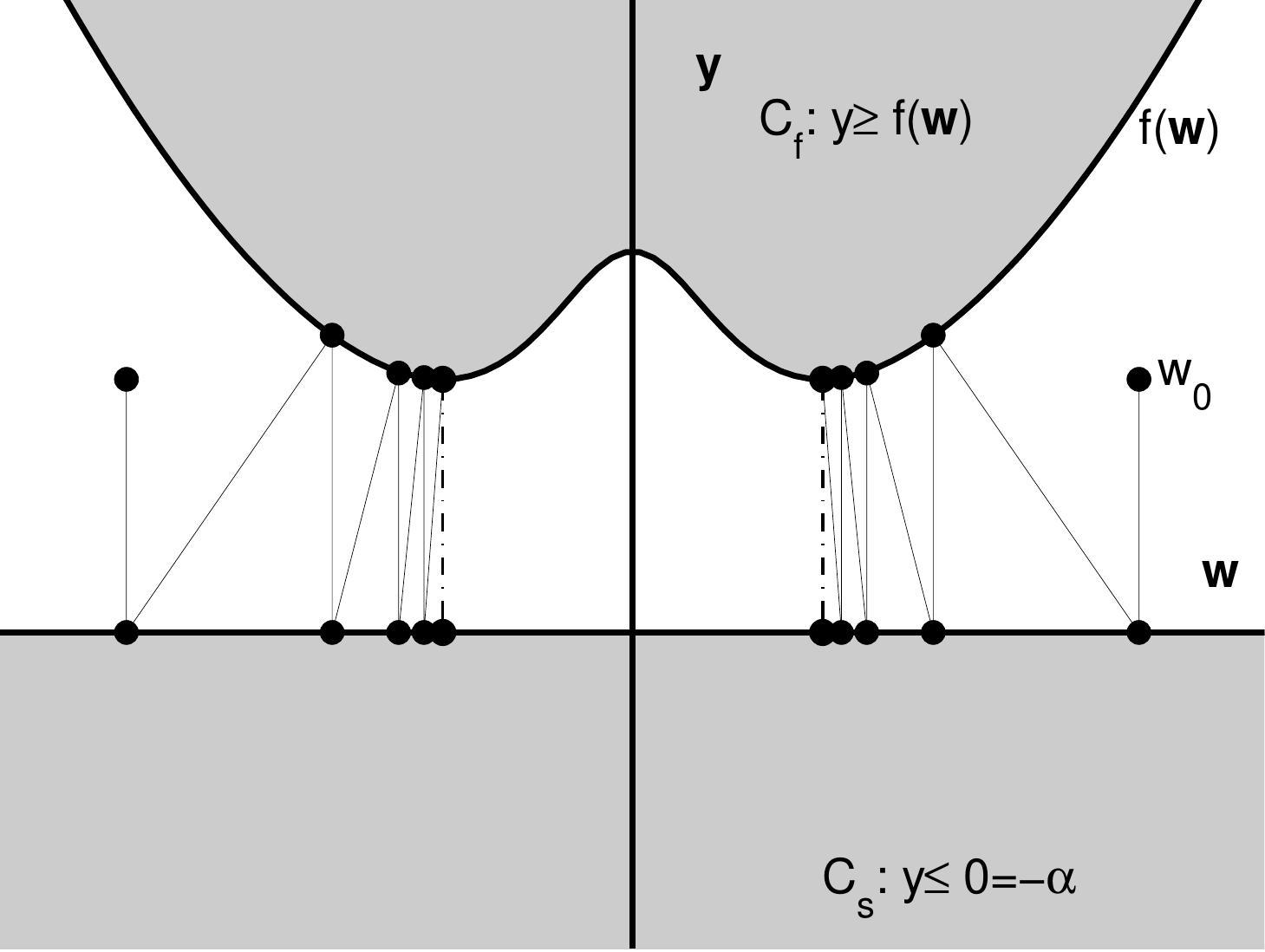}
\caption[Converging to the cost function.]{Iterations start with an initial vector $\mathbf{w_0}$ and iterates converge to a local minimum by the POCS algorithm.}
\label{app:noconvex}
\end{figure}

If the cost function $f$ is not convex and have more than one local minimum then the corresponding set $C_f$ is not convex in $R^{N+1}$. In this case iterates may converge to one of the local minima. This is graphically illustrated in Fig. \ref{app:noconvex}.

\section{Supporting Hyperplane Concept based POCS Solution}
It may not be easy to find the orthogonal projection onto the set C$_{\mathrm{f}}$ for some cost functions $f$. In such cases it is possible to use supporting hyperplanes of the convex set to find the minimum of $f(\mathbf{w})$. The second optimization algorithm is based on making successive orthogonal projections onto the supporting hyperplanes of the set $C_f$ instead of the actual set. 

\begin{figure}[t!]
\includegraphics[width=\linewidth,clip,trim=0.25in 0.25 0.25 0.25in]{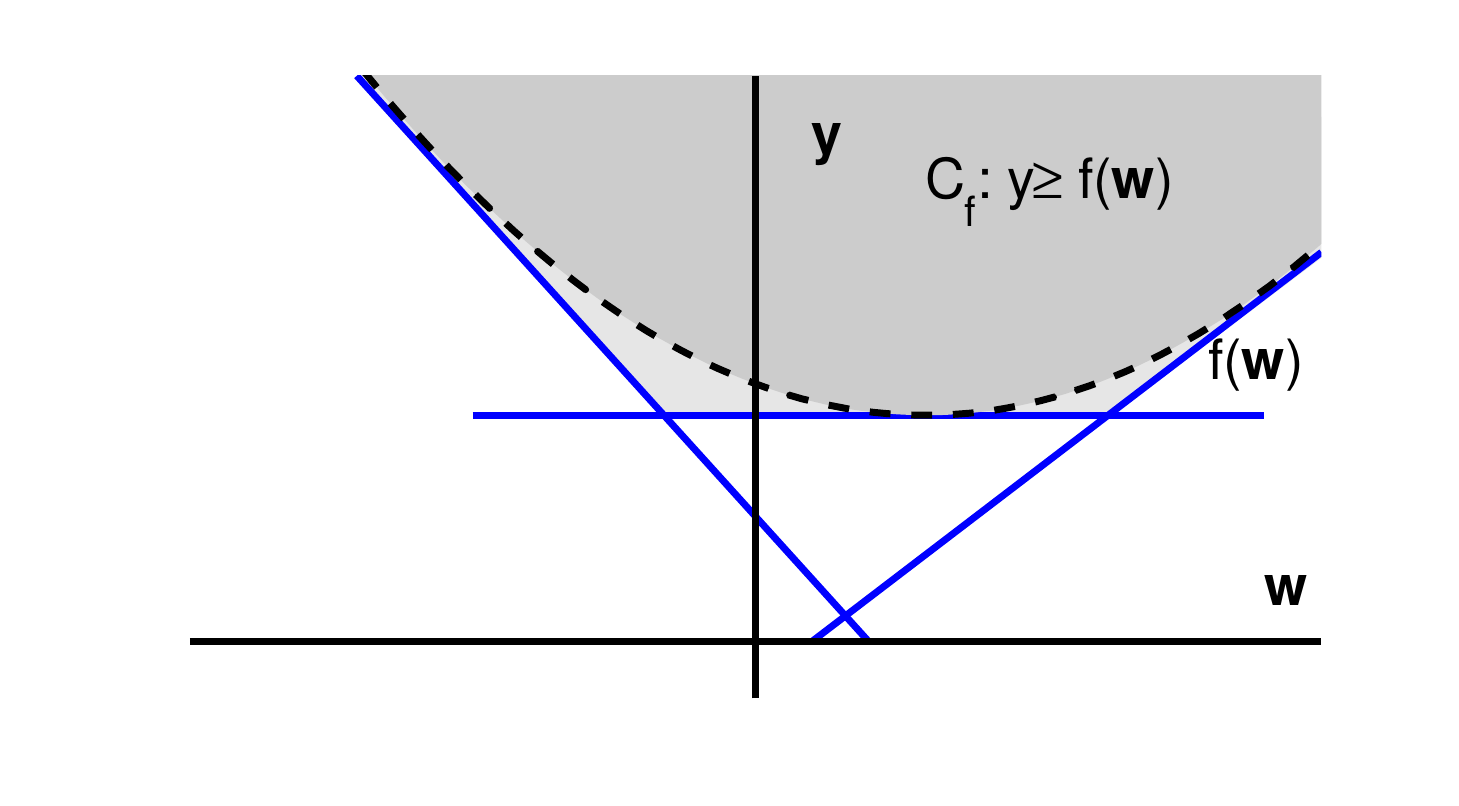}
\caption[Intersection of half planes.]{Intersection of supporting halfplanes C$_{\mathrm{f},\mathbf{w}}$ defines the set C$_{\mathrm{f}} = \bigcap\limits_{\mathbf{w} \in \mathbb{R}} \mathrm{C_{f,\mathbf{w}}}$ when the function f is convex.}
\label{app:supporthyp}
\end{figure}

The set C$_{\mathrm{f}}$ can be expressed as the intersection of halfplanes whose boundaries are supporting hyperplanes as shown in Fig. \ref{app:supporthyp}.
Let $\mathbf{w}$ and $f(\mathbf{w})$ form a vector $ [ \mathbf{w}^T f(\mathbf{w}) ]^T $ in $\mathbb{R}^{N+1}$ on the surface of $y=f(\mathbf{w})$. Let the supporting hyperplane at this point be $l(\mathbf{w})$. Let us also define the halfplane (or halfspace) set as follows:
\begin{equation}
\label{app:eq:c8}
\mathrm{C_{f,\mathbf{w}}} = \{y \geq l(\mathbf{w})\}
\end{equation}
Clearly, the set C$_{\mathrm{f}}$ can be expressed as the intersection of its supporting halfspaces in $R^{N+1}$: 
\begin{equation}
\label{app:eq:c9}
\mathrm{C_f} = \bigcap\limits_{\mathbf{w} \in \mathbb{R}^N} \mathrm{C_{f,\mathbf{w}}}
\end{equation}

Therefore, the POCS approach can be applied to the level set C$_{\mathrm{s}}$ and the family of sets C$_{\mathrm{f},\mathbf{w}}$, $\mathbf{w} \in \mathbb{R}^N$ to find the minimum of $f(\mathbf{w})$. In this case, the number of sets are infinite. This set theoretic scenario was studied by Slavakis, Yukawa, Yamada and Theodoridis \cite{Sla08,Sla09}.

\begin{figure}[t!]
\includegraphics[width=\linewidth,clip,trim=0.25in 0.25 0.25 0.25in]{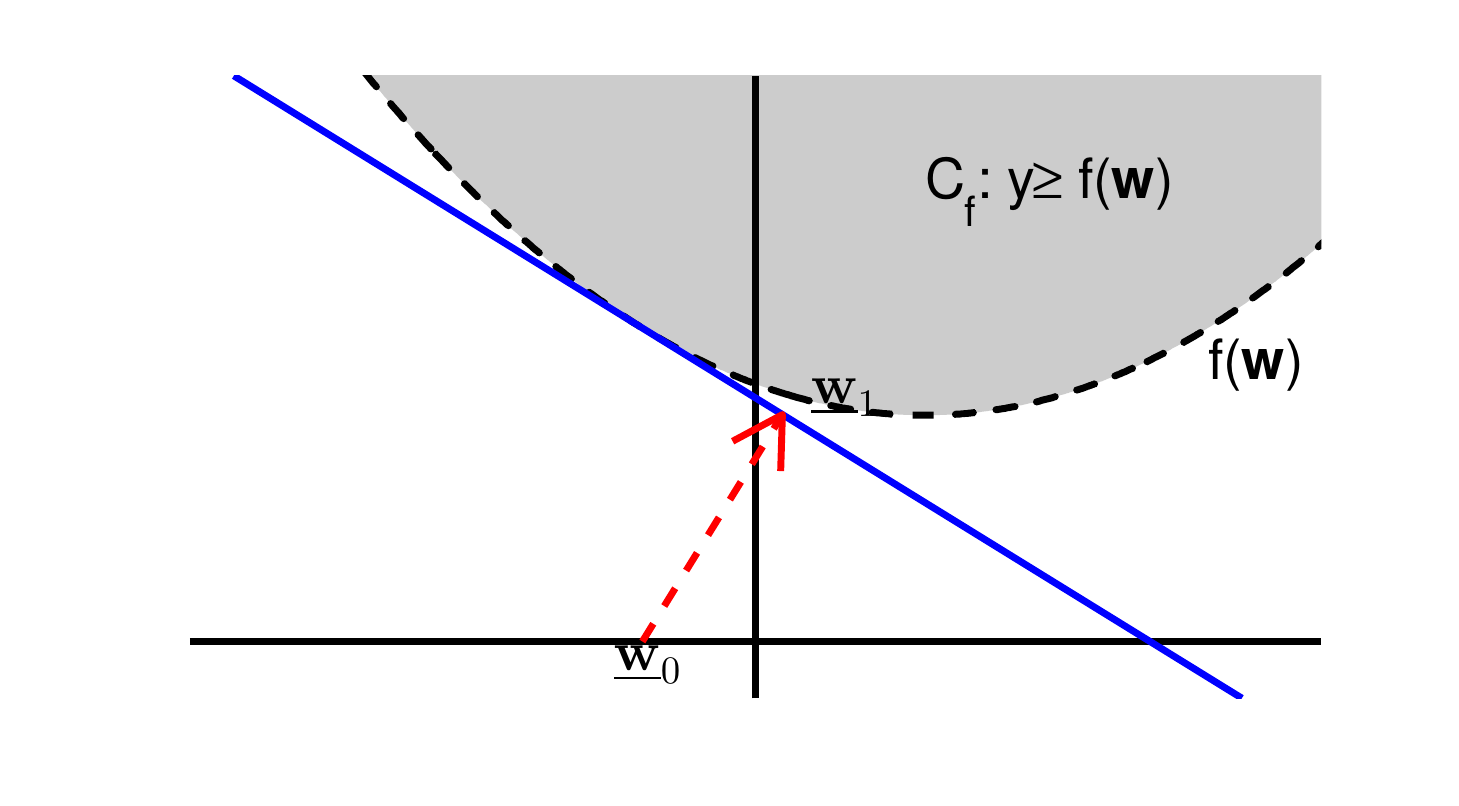}
\caption[Projection $\underline{w}_1$ of the vector $\underline{w}_o$ onto a supporting hyperplane.]{Projection $\underline{w}_1$ of the vector $\underline{w}_o$ onto a supporting hyperplane}.
\label{app:projsupporthyp}
\end{figure}

In the second optimization approach we perform orthogonal projections onto supporting hyperplanes of the cost function instead of the actual set $\mathrm{C_f}$ as shown in Fig. \ref{app:projsupporthyp}. 
Since making an orthogonal projection onto a hyperplane is easy to compute, the optimization problem does not require the solution of any nonlinear equations as long as it is possible to compute the surface normal at $f(\mathbf{w})$. Let the surface normal at $[ \bf{w}_t^T ~  f(\bf{w}_t)]^T$ be
$ \bf{v}_t \in R^{N+1}$. The supporting hyperplane $l(\bf{\underline{w}}_t)$ is given by 
\begin{equation}
\label{app:support}
\mathbf{\underline{v}}_t^T \dot ({\underline{\bf{w}}} - {\underline{\bf{w}}}_t) =0
\end{equation}
The orthogonal projection ${\underline{\bf{w}}}_p$ of an arbitary vector 
${\underline{\bf{w}}}_o$ onto the hyperplane  $l(\bf{\underline{w}}_t)$ is obtained as follows
\begin{equation}
\label{app:projection}
\mathbf{\underline{w}}_p =  {\underline{\bf{w}}_o} ~ + ~ \lambda ~ \frac{\mathbf{\underline{v}}_t^T \dot ({\underline{\bf{w}}_o} - {\underline{\bf{w}}}_t) }{||\mathbf{\underline{v}}_t||^2} \mathbf{\underline{v}}_t
\end{equation}
where $\lambda =1$. The parameter  $\lambda$ can be selected between  
 $0< \lambda <2$ as in the normalized LMS algorithm \cite{Sla08}.  The parameter  $\lambda \neq 1$ case corresponds to non-orthogonal projections.
Eq. (\ref{app:projection}) is the key equation of the supporting hyperplane based optimization approach. In Fig. \ref{x} a graphical illustration of the iterative optimization algorithm is shown.
\begin{figure}[ht!]
\label{x}
\centering
\includegraphics[width=0.8\linewidth]{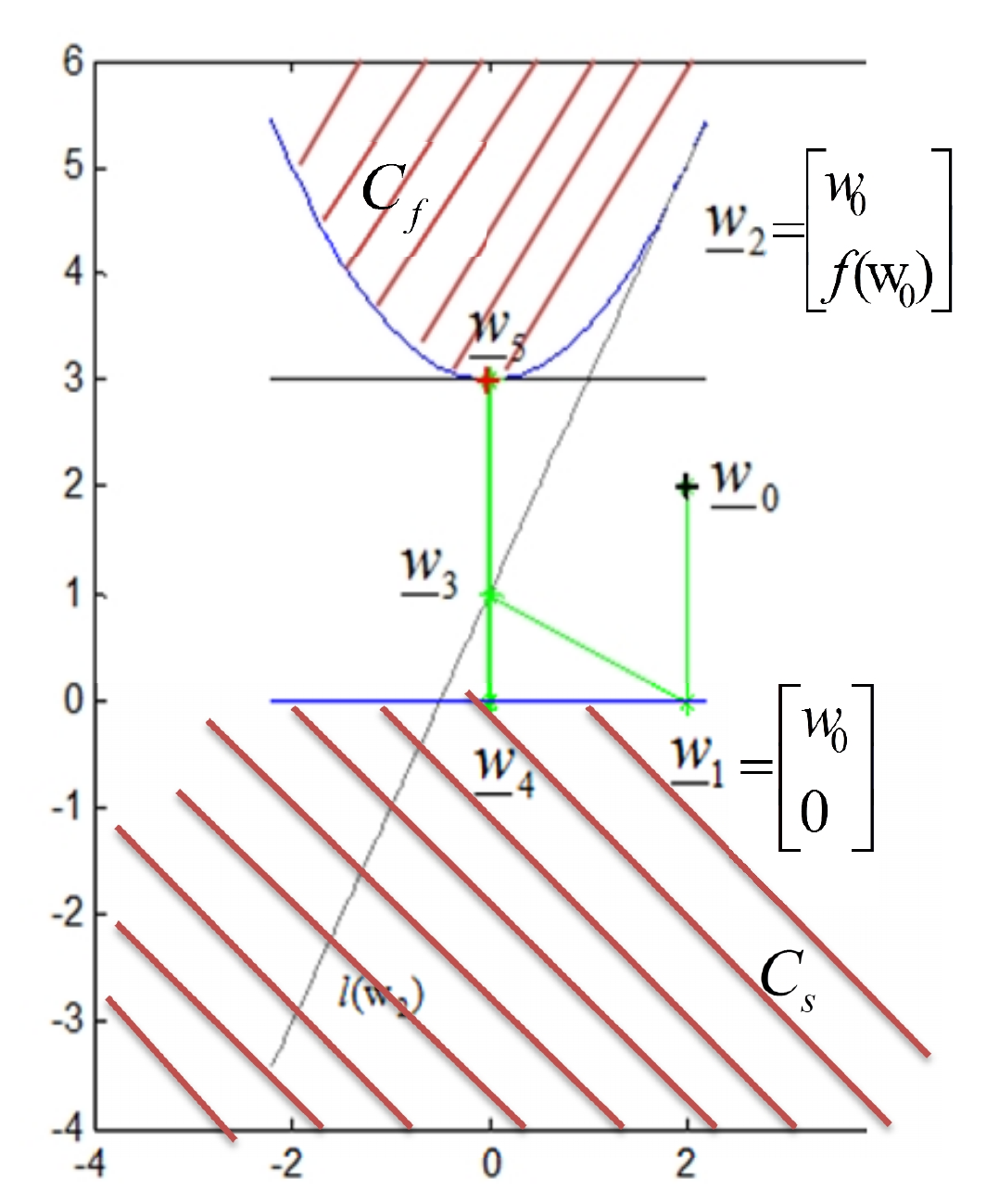}
\caption[Projections onto supporting hyperplanes of the set $C_f$ lead to the global minimum.]{Projections onto supporting hyperplanes of the set $C_f$ lead to the global minimum}.
\label{app:pocs21}
\end{figure}
 Iterations start with an arbitrary $\underline{\bf{w}}_o \in R^{N+1} $. The vector $\underline{\bf{w}}_o$ 
is projected onto the set $C_s$ and $\underline{\bf{w}}_1$  is obtained. The projected vector is 
\begin{equation}
\underline{\bf{w}}_1 = [ {\bf w}_o^T ~ ~ 0 ]^T
\end{equation}
where $\bf{w}_o$ is an $N$ dimensional vector containing the first $N$ components of the $N+1$ dimensional vector $\underline{\bf{w}}_o$. 
The value of the cost function $f(\bf{w}_o)$ and its surface normal 
 $\underline{\bf{v}}_o$ at  $\bf{w}_o$ are computed. This is the next iterate
\begin{equation}
\underline{\bf{w}}_2 = [  {\bf w}_o^T ~ ~ f ({\bf w}_o) ]^T
\end{equation}
The corresponding supporting hyperplane $l(\underline{\bf{w}}_2)$ is characterized by the equation:
\begin{equation}
\label{app:support}
\mathbf{\underline{v}}_o^T \dot ({\underline{\bf{w}}} - {\underline{\bf{w}}}_o) =0
\end{equation}
The next iterate $\underline{\bf{w}}_3$ is determined by projecting 
$\underline{\bf{w}}_1$ onto the hyperplane $l(\underline{\bf{w}}_2)$ as follows:
\begin{equation}
\label{app:projection2}
\mathbf{\underline{w}}_3 =  {\underline{\bf{w}_1}} + ~ \lambda ~ \frac{\mathbf{\underline{v}}_t^T \dot ({\underline{\bf{w}}_1} - {\underline{\bf{w}}}_t) }{||\mathbf{\underline{v}}_o||^2} \mathbf{\underline{v}}_o
\end{equation}

\begin{figure}[ht!]
\centering
\includegraphics[width=0.7\linewidth,clip,trim=0.25in 0.25 0.25 0.25in]{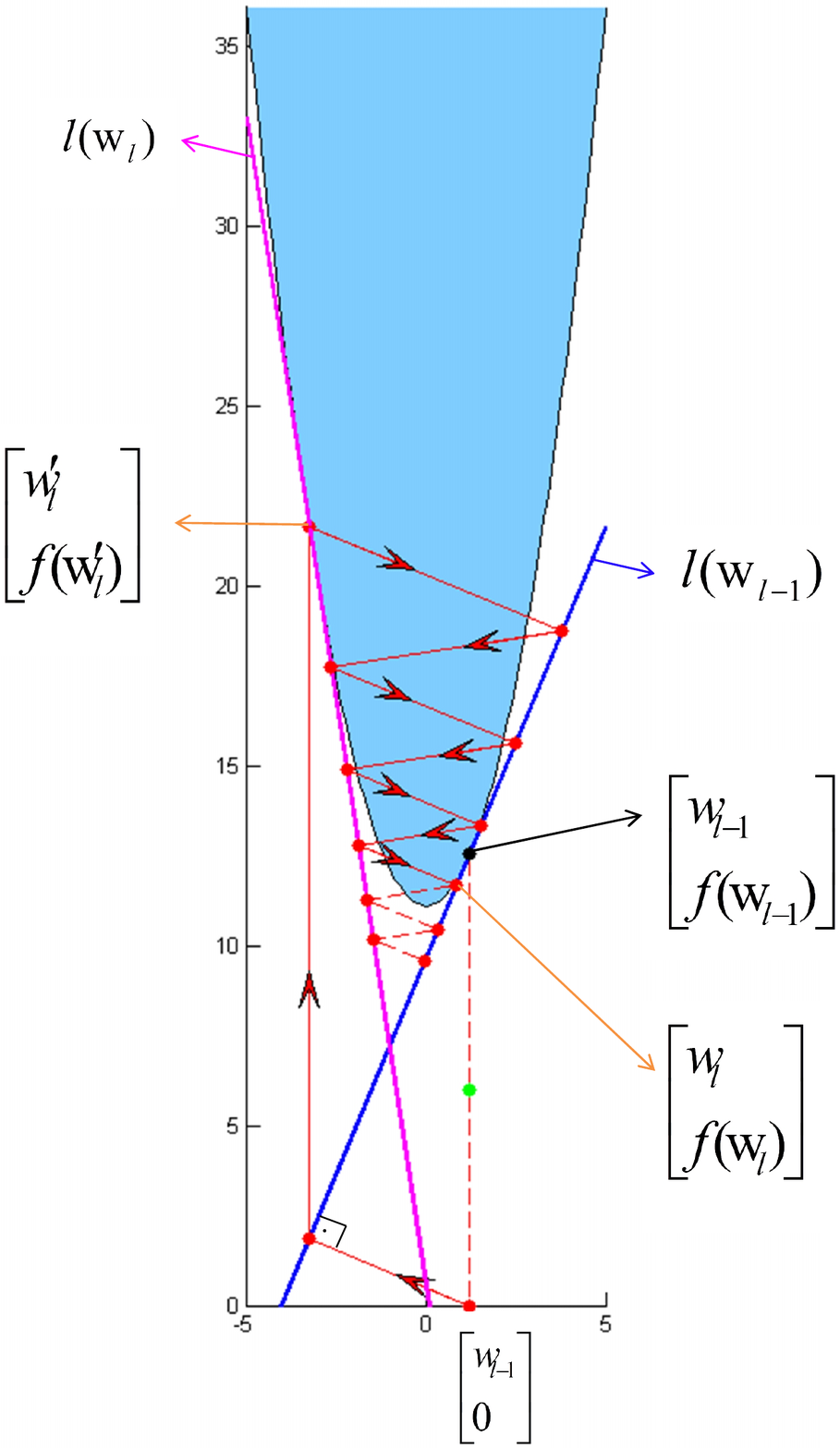}
\caption[]{Projections onto two supporting hyperplanes $l(\bf w'_l)$ and  $l(\bf w_{l-1})$ to obtain the next iterate satisfying $f(\bf w_l) > f(\bf w_{l-1})$.}
\label{app:twoHPs}
\end{figure}

This computes the first iteration cycle. Next, we project $\mathbf{\underline{w}}_3$ onto the set $C_s$ and 
$\underline{\bf{w}}_4 = [  {\bf w}_3^T ~ ~ 0 ]^T$ is obtained. The vector  $\bf{w}_3$ is an $N$ dimensional vector containing the first $N$ components of the $N+1$ dimensional vector $\underline{\bf{w}}_4$. We obtain 
$\underline{\bf{w}}_5 = [  {\bf w}_3^T ~ ~ f({\bf w}_3) ]^T$ on the surface of the set $C_f$.
At this point, we verify if $f(\bf w_3)$ is less than $f(\bf w_o)$ or not. If yes, we continue the iterations as described above.
If not, we switch to another iteration strategy as graphically shown in Fig. \ref{app:twoHPs}. 

Let us assume that $f(\bf w'_l) > f(\bf w_{l-1})$ as shown in
Fig.  \ref{app:twoHPs}. In this case consecutive iterative projections are performed onto the supporting hyperplanes  $l(\bf w'_l)$ and  $l(\bf w_{l-1})$
until a vector $[  {\bf w}_l^T ~ ~ f({\bf w}_l) ]^T$ satisfying
$f(\bf w_l) > f(\bf w_{l-1})$ is obtained.


Projection onto a supporting hyperplane is easy to perform. However the cost function may not have a well defined derivative at a given vector ${\bf w}_ud$ and a well defined supporting hyperplane may not exist. In this case any hyperplane $l({\bf w_d})$ passing through ${\bf w}_ud$ and satisfying $f({\bf w}) ~ \ge ~l({\bf w_d})$ can be used in Eq. \ref{app:projection2}. 


It is also possible to include other convex constraints into the optimization problem:
\begin{eqnarray}
min ~ ~ f ({\bf w}) \\
such ~ that~  {\bf w} \in C_1, C_2 , \ldots , C_L \nonumber
\label{constraint}
\end{eqnarray}
where $C_1, C_2 , \ldots , C_L$ are closed and convex sets representing constraints on the solution of the inverse problem. However 
 the main difficulty with this approach is that non-intersecting multiple convex set scenario has not been fully studied to the best of our knowledge. Successive orthogonal projections onto non-intersecting convex sets may lead to limit cycles \cite{Gub67,Cen12}. This remains as an interesting research problem. It is experimentally observed that successive orthogonal projections onto $C_f$, $C_s$ 
and another constraint set leads to a limit cycle containing the optimal solution. 

One possible way to handle the problem described in (\ref{constraint}) is to enlarge one or some of the sets so that they have a well defined non-empty intersection. We successfully applied this strategy in FIR filter design \cite{Cet97}. For example, it is trivial to enlarge the set $C_s$ by slowly increasing the value of $\alpha$ in a judicuous manner. 




It is possible to develop  iteration strategies based on non-orthogonal projections or linear combinations of projections to speed up the convergence to the global minimum as discussed in \cite{Com12,Com93,Com04}.


\section{Minimization of Non-Convex Functions}

\begin{figure}[ht!]
\label{x}
\centering
\includegraphics[width=0.8\linewidth]{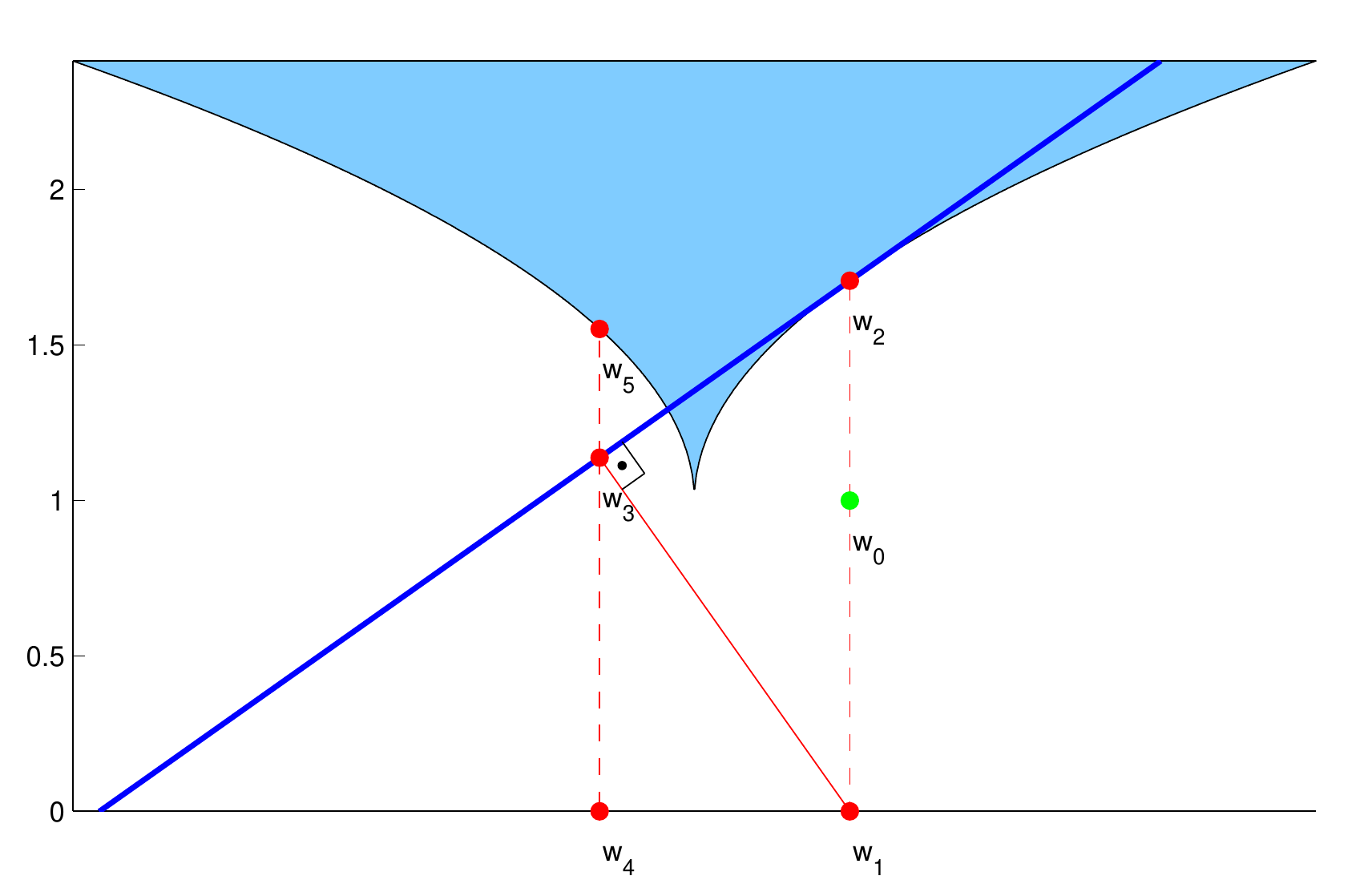}
\label{app:pocs21}
\caption{ The function $f(w) = |w-2|^{0.5} + 1$ is not convex. However projections onto tangential hyperplanes may lead to the global minimum.}
\end{figure}

An important class of cost functions are based on  $l_p,$ $p<1$. It is experimentally observed that such functions can be minimized by using the "supporting" hyperplane concept as shown in Fig. \ref{app:pocs21}.  Obviously, tangential hyperplanes are no longer "supporting" hyperplanes and the set $C_s$ is not a convex set. This is the "inside out" version of the convex minimization problem that is studied in Section 3 and the iterative scheme introduced in Section 3 leads to the minimum of the function. 

P. Combettes provided an excellent review of POCS theory and the recently introduced proximal splitting method in \cite{Com12,Com11}. The relation between the proposed methods and the proximal splitting theory will be investigated in the future.

Another interesting future research direction is the use of generalized
 supporting hyperplane concept to minimize non-convex functions with many local minima. The convex hull of the generalized supporting planes may form a convex region in $\mathbb{R}^{N+1}$. As a result it may be possible to find the minimum of the cost function by performing successive orthogonal projections onto generalized supporting hyperplanes. This problem will be also studied in the future.




\bibliographystyle{IEEEtranIO}
\bibliography{PhdReferences}

\end{document}